\documentclass[11pt,a4paper]{article}

\usepackage{graphicx}
\usepackage{latexsym}
\usepackage{amsmath, amsfonts,amssymb, amsthm, euscript,makeidx, xcolor,mathrsfs,algorithm}
\usepackage{enumerate}
\usepackage{tabu}

\usepackage[colorlinks,linkcolor=blue,anchorcolor=green,citecolor=red]{hyperref}

\def\sqr#1#2{{\vcenter{\vbox{\hrule height.#2pt
              \hbox{\vrule width.#2pt height#1pt \kern#1pt \vrule width.#2pt}
              \hrule height.#2pt}}}}
\def\5n{\negthinspace \negthinspace \negthinspace \negthinspace \negthinspace }
\def\4n{\negthinspace \negthinspace \negthinspace \negthinspace }
\def\3n{\negthinspace \negthinspace \negthinspace }
\def\2n{\negthinspace \negthinspace }
\def\1n{\negthinspace }

\def\Ker{\operatorname{Ker}}
\def\Ran{\operatorname{Ran}}
\def\Rank{\operatorname{Rank}}

\def\dbF{\mathbb{F}}

\def\dbN{\mathbb{N}}

\def\dbP{\mathbb{P}}

\def\dbR{\mathbb{R}}
\def\dbS{\mathbb{S}}

\def\={\buildrel \triangle \over =}

\def\b{\beta}

\def\e{\varepsilon}

\def\l{\lambda}

\def\f{\varphi}
\def\th{\theta}

\def\mf{\mathcal{F}}
\def\me{\mathbb{E}}
\def\rd{\,\mathrm d}
\def\bal{\begin{aligned}}
\def\eal{\end{aligned}}
\def\cA{{\cal A}}

\def\cC{{\cal C}}

\def\cF{{\cal F}}

\def\cJ{{\cal J}}
\def\cK{{\cal K}}
\def\cL{{\cal L}}
\def\cM{{\cal M}}
\def\cN{{\cal N}}

\def\cT{{\cal T}}
\def\cU{{\cal U}}
\def\cV{{\cal V}}

\def\BA{{\bf A}}

\def\BK{{\bf K}}
\def\BL{{\bf L}}

\def\Bv{{\bf v}}

\def\Bx{{\bf x}}

\def\Bz{{\bf z}}
\def\no{\noindent}

\def\ms{\medskip}

\def\q{\quad}
\def\qq{\qquad}

\def\lan{\langle}
\def\ran{\rangle}

\def\rf{\eqref}

\def\h{\widehat}
\def\wt{\widetilde}

\def\cd{\cdot}
\def\cds{\cdots}

\def\ae{\hbox{\rm a.e.}}
\def\as{\hbox{\rm a.s.}}

\def\deq{\triangleq}

\def\({\Big (}
\def\){\Big )}

\def\bde{\begin{definition}\label}
\def\ede{\end{definition}}
\def\be{\begin{equation}}
\def\bel{\begin{equation}\label}
\def\ee{\end{equation}}
\def\beq{\begin{equation*}\begin{aligned}}
\def\eeq{\end{aligned}\end{equation*}}
\def\bt{\begin{theorem}\label}
\def\et{\end{theorem}}
\def\bc{\begin{corollary}\label}
\def\ec{\end{corollary}}
\def\bl{\begin{lemma}\label}
\def\el{\end{lemma}}
\def\bp{\begin{proposition}\label}
\def\ep{\end{proposition}}
\def\bas{\begin{assumption}\label}
\def\eas{\end{assumption}}
\def\br{\begin{remark}\label}
\def\er{\end{remark}}
\def\bex{\begin{example}\label}
\def\ex{\end{example}}
\def\ba{\begin{array}}
\def\ea{\end{array}}
\def\ed{\end{document}}

\def\square#1{\vbox{\hrule\hbox{\vrule height#1%
     \kern#1\vrule}\hrule}}
\def\rectangle#1#2{\vbox{\hrule\hbox{\vrule height#1%
     \kern#2\vrule}\hrule}}

\font\tenbb=msbm10 \font\sevenbb=msbm7 \font\fivebb=msbm5

\newfam\bbfam
\scriptscriptfont\bbfam=\fivebb \textfont\bbfam=\tenbb
\scriptfont\bbfam=\sevenbb

\allowdisplaybreaks  

\newtheorem{theorem}{\hskip 1.3em Theorem}[section]
\newtheorem{definition}[theorem]{\hskip 1.3em Definition}
\newtheorem{proposition}[theorem]{\hskip 1.3em Proposition}
\newtheorem{corollary}[theorem]{\hskip 1.3em Corollary}
\newtheorem{lemma}[theorem]{\hskip 1.3em Lemma}
\newtheorem{remark}[theorem]{\hskip 1.3em Remark}
\newtheorem{example}[theorem]{\hskip 1.3em Example}

\newtheorem{assumption}[theorem]{\hskip 1.3em Assumption}

\makeatletter
   
   \@addtoreset{equation}{section}
\makeatother

\begin{document}

\title{Minimal norm control problems for linear mean-field control systems}

\author{Manfen An,\thanks{
School of Mathematics and Statistics, Southwest University, Chongqing 400715, China. {\small\it
e-mail:} {\small\tt xndx135956@email.swu.edu.cn}. \ms}
~~~Rui Gao,\thanks{School of Mathematics and Statistics, Southwest University, Chongqing 400715, China.  {\small\it
e-mail:} {\small\tt rogerui@swu.edu.cn}. \ms}~~~
  and~~~
Yanqing Wang,\thanks{Corresponding author,
School of Mathematics and Statistics, Southwest University, Chongqing 400715, China.  {\small\it
e-mail:} {\small\tt yqwang@swu.edu.cn}.}}
\date{}
\maketitle

\begin{abstract}
This paper studies controllability and minimal norm control problems (NP problems) for mean-field stochastic differential equation (MF-SDE) systems. We first establish an equivalent rank criterion for terminal controllability, namely, that the diffusion matrices $D$ and $D+\bar D$ are both of full row rank. 
We then derive an explicit Riccati-based representation of the minimal norm control by applying a parameterization technique. A key step is to identify the value $P(0)$ of an associated Riccati solution with the controllability Gramian of a controllable auxiliary system, which yields the parameter required in the construction. The resulting representation provides a basis for the future development and analysis of numerical methods for NP problems.
\end{abstract}

\ms

\no\bf Keywords: \rm 
mean-field control system,
controllability,
minimal norm control problem
\ms

\no\bf AMS 2020 subject classification: \rm 
60H10,
93B05,
93E20

\section{Introduction}

This paper is concerned with the minimal norm control problem (NP problem, for short) for the following linear mean-field stochastic differential equation (MF-SDE, for short):
\bel{fmf}
\bal
& \mathrm dx(t)=\{Ax(t)+\bar A\me[x(t)]
+Bu(t)+\bar B\me[u(t)]\}\,\mathrm dt \\
&\qquad\q+
\{Cx(t)+\bar C\me[x(t)]
+D u(t)+\bar D\me[u(t)]\}\,\mathrm dW(t),\q t\in[0,T], \\
\eal
\ee
where $T\in (0,\infty)$, $W(\cd)$ is a one-dimensional standard Brownian motion, $A,\bar A, C,\bar C \in\dbR^{n\times n}$, $B,\bar B, D,\bar D \in\dbR^{n\times m}$ with $n,m\in\dbN$, and $u(\cd)\in L^2_\dbF(0,T;\dbR^m)$. 
For any $x_0\in\dbR^n$ and any $x_T\in L^2_{\mf_T}(\Omega;\dbR^n)$, we consider the following NP problem.\\
{\bf Problem (NP)$_{x_0,x_T}^T$.} Find a control $u^*(\cd)\in \cU(T,x_0,x_T)$ minimizing the cost functional
\bel{cost-F}
\cJ_F\big(u(\cd)\big)=\|u(\cd)\|^2_{L^2_\dbF(0,T;\dbR^m)},
\ee
where the admissible control set is defined by
\bel{ad-control}
\cU(T,x_0,x_T)
\deq \big\{
u(\cdot)\in L^2_{\dbF}(0,T;\dbR^m)
\,\big|\,
x\big(T;x_0,u(\cdot)\big)=x_T,\quad \mbox{$\dbP$-a.s.}
\big\}.
\ee
Here $x\big(\cd;x_0,u(\cdot)\big)$ denotes the solution of MF-SDE \rf{fmf} with the initial condition $x(0)=x_0$ and the control $u(\cd)$. 
The spaces $L^2_{\mf_T}(\Omega;\dbR^n)$ and $L^2_\dbF(0,T;\dbR^m)$ are defined in Section \ref{subsec-notation}.

To study Problem {\bf (NP)$_{x_0,x_T}^T$} for arbitrary initial-terminal pairs $(x_0,x_T)$, one first needs to guarantee that the admissible control set $\cU(T,x_0,x_T)$ is nonempty for every such pair. This is precisely the exact controllability of MF-SDE \rf{fmf}; see Definition \ref{def-control} in Section \ref{subsec-NP}. Therefore, exact controllability provides the fundamental solvability condition for the minimal norm control problem considered in this paper.

For stochastic differential equation (SDE) systems, which correspond to the special case $\bar A=\bar C=0$ and $\bar B=\bar D=0$, the study of controllability in a probabilistic sense can be traced back to \cite{Sunahara-Aihara-Kishino75,Zabczyk81}.
With the aid of backward stochastic differential equations,
 {\em exact controllability} and {\em terminal controllability} in the $L^2$ sense were subsequently studied in \cite{Peng94}, where corresponding rank conditions were also derived. In particular, a necessary condition for terminal controllability is 
\bel{wang925e4}
\Rank D=n.
\ee
Later, exact and partial controllability in the $L^p$ sense for SDE systems with time-varying coefficients were investigated in, e.g.,
\cite{Liu-Peng10,Lv-Yong-Zhang12, Wang-Yang-Yong-Yu17,Wang-Yu20}.
On the other hand, {\em approximate controllability} in the $L^2$ sense was studied in
\cite{Buckdahn-Quincampoix-Tessitore06,Goreac08}.

For the MF-SDE system \rf{fmf}, exact controllability has been studied in, e.g., \cite{Goreac14,Yu21,Liu-Yu23, Chen-Yu24,Goreac-Li-Zhang26}.
The full-row-rank condition
\bel{wang925e5}
\Rank D=\Rank(D+\bar D)=n
\ee
extends the corresponding rank condition \rf{wang925e4} for SDE systems. In the existing literature, condition \rf{wang925e5} appears either as a structural assumption (see
\cite[Proposition 5]{Goreac14},
\cite[Assumption (H)]{Yu21}, 
\cite[Assumption 7]{Liu-Yu23}, 
and
\cite[Assumption (H4$'$)]{Chen-Yu24}), or is derived under additional restrictions on the coefficients, such as $C=0$ in \cite{Goreac-Li-Zhang26}.
In the present work, we close this gap by proving, without imposing any additional restrictions on the coefficients, that
\beq
\mbox{MF-SDE \rf{fmf} is terminally controllable}
\quad\Longleftrightarrow\quad
\Rank D=\Rank(D+\bar D)=n;
\eeq
see Theorem \ref{rank-condition}. This constitutes the first main contribution of this paper.

The above controllability analysis addresses the feasibility of Problem {\bf (NP)$_{x_0,x_T}^T$}. We next turn to the optimization issue: once the admissible control set is nonempty, how can the minimal norm control $u^*(\cd)$ of Problem {\bf (NP)$_{x_0,x_T}^T$} be characterized and explicitly constructed?

The NP problems considered in this paper belong to the class of linear-quadratic optimal control problems (LQ problems, for short) with state constraints. NP problems for SDE systems were first studied in \cite{Wang-Zhang15,Gashi15}, and were further developed in \cite{Wang-Yang-Yong-Yu17,Bi-Sun-Xiong20,Wang-Yu20,Wang26}. In particular, the equivalence between the exact controllability of the underlying SDE systems and the unique solvability of the corresponding NP problems was established; see \cite[Theorems 4.1, 4.5, 5.1]{Wang-Yang-Yong-Yu17} and \cite[Theorem 2.5]{Wang-Yu20}.

For the MF-SDE system \rf{fmf}, exact controllability guarantees the nonemptiness of the admissible control set $\cU(T,x_0,x_T)$ for every prescribed initial-terminal pair. Since $\cU(T,x_0,x_T)$ is a closed affine subset of the Hilbert space $L^2_\dbF(0,T;\dbR^m)$, the corresponding NP problem then admits a unique minimal norm control. The main issue is therefore not merely the existence and uniqueness of the minimal norm control, but rather how to characterize and construct it explicitly.
In the present work, we employ a parameterization approach that yields an explicit Riccati-based representation of the minimal norm control. The main steps are outlined as follows.
\begin{enumerate}
\item [\rm (i)] Under the exact controllability assumption on \rf{fmf}, we first apply the rank condition \rf{wang925e5} established above to transform the original MF-SDE system \rf{fmf} into the equivalent mean-field backward stochastic differential equation (MF-BSDE) system \rf{fmf-2}. We then show that the solvability of Problem {\bf (NP)$_{x_0,x_T}^T$} for \rf{fmf} is equivalent to that of Problem {\bf (COP)$_{x_0,x_T}^T$} for \rf{fmf-2}; see Section \ref{subsec-transformation}.

\item [\rm (ii)] Since Problem {\bf (COP)$_{x_0,x_T}^T$} is subject to the constraint $x(0)=x_0$, we introduce a parameterized unconstrained LQ problem, Problem {\bf (OP)$_{x_0,x_T,\l}^T$}, associated with \rf{fmf-2}. We show that, if there exists a parameter $\l\in\dbR^n$ such that
\bel{wang925e8}
x^*_\l(0)=x_0,
\ee
then the corresponding optimal control of Problem {\bf (OP)$_{x_0,x_T,\l}^T$} is the unique optimal control of Problem {\bf (COP)$_{x_0,x_T}^T$}; see Section \ref{subsec-para}.

\item [\rm (iii)]
We apply the Riccati decoupling technique to Problem {\bf (OP)$_{x_0,x_T,\l}^T$} and derive feedback representations for its optimal state and control. By constructing an auxiliary controllable system and identifying its controllability Gramian with the value $P(0)$ of the associated Riccati solution, we prove the invertibility of $P(0)$ and thereby determine $\l$ satisfying \rf{wang925e8}. This ultimately yields an explicit Riccati-based representation of the minimal norm control of Problem {\bf (NP)$_{x_0,x_T}^T$}; see Section \ref{subsec-para-choice}.
\end{enumerate}
The resulting representation provides a constructive characterization of the minimal norm control in terms of Riccati equations, an MF-BSDE, and a decoupled MF-SDE. It also provides a basis for the development and analysis of numerical algorithms for NP problems. This constitutes the second main contribution of this paper.


The remainder of this paper is organized as follows. Section \ref{sec:prelim} introduces the notation, recalls the notions of terminal and exact controllability for MF-SDE systems, and defines minimal norm controls associated with NP problems. In Section \ref{sec-rank}, we establish the equivalent rank criterion for terminal controllability. Section \ref{sec-representation} is devoted to the Riccati-based representation of the minimal norm control. We first transform the original problem into an equivalent constrained optimal control problem for an MF-BSDE, then introduce a parameterized unconstrained LQ problem, and finally determine the parameter through a controllability Gramian and derive the explicit representation of the minimal norm control.

\section{Preliminaries}\label{sec:prelim}
	
\subsection{Notation}\label{subsec-notation}

Let $\dbN$ be the set of positive integers.
For $n,m\in\dbN$, denote by $\dbR^n$ the $n$-dimensional Euclidean space equipped with the standard inner product $\lan\cd,\cd\ran$ and the induced norm $\|\cd\|$, 
by $\dbR^{n\times m}$ the space of all $n\times m$ real matrices,
by $\dbS^n$ the space of all $n\times n$ real symmetric matrices, 
and by $\dbS^n_+$ the cone of all positive semidefinite matrices in $\dbS^n$.
Throughout this paper, the superscript ``$\top$'' denotes the transpose of a vector or a matrix.

Let $T>0$ be a fixed terminal time. Let $(\Omega,\cF,\dbF,\dbP)$ (with $\dbF=\{\cF_t\}_{t\in[0,T]}$) be a complete filtered probability space, on which a one-dimensional standard Brownian motion $\{W(t)\}_{t\in[0,T]}$ is defined such that $\dbF$ is the natural filtration generated by $W(\cd)$ (augmented by all $\dbP$-null sets).

Throughout this paper, we use the following spaces:
\begin{align*}
&L^2_{\cF_T}(\Omega;\dbR^n)=\big\{\xi:\Omega\to \dbR^n\,\big|\, \xi \mbox{ is $\mf_T$-measurable, and } \me\big[\|\xi\|^2\big]<\infty\big\},\\
&L^2_\dbF(0,T;\dbR^n)=\Big\{x:[0,T]\times \Omega\to \dbR^n\,\big|\, x(\cd) \mbox{ is $\dbF$-adapted, and } \me\!\int_0^T\!\|x(t)\|^2\rd t<\infty\Big\},\\
&L^2_\dbF(\Omega;C([0,T];\dbR^n))=\Big\{x(\cd)\in L^2_\dbF(0,T;\dbR^n)\,\big|\, x(\cd) \mbox{ has continuous paths,}\\
&\qq\qq\qq\qq\qq\qq \mbox{ and } \me\Big[\sup_{t\in[0,T]}\|x(t)\|^2\Big]<\infty\Big\}.
\end{align*}


For any square-integrable random variable or process $\xi$, we denote its mean by $\bar\xi$ and its centered counterpart by $\widetilde\xi$, namely,
\beq
\bar\xi=\me[\xi],
\qquad
\widetilde\xi=\xi-\me[\xi].
\eeq
Thus, $\xi=\widetilde\xi+\bar\xi$ and $\me[\widetilde\xi]=0$.
For any two matrices $H$ and $\bar H$ of the same size, we introduce the notation
\beq
\widehat H=H+\bar H.
\eeq
	
\subsection{Exact controllability and minimal norm problem for mean-field systems}\label{subsec-NP}
	
For every $u(\cd)\in L^2_{\dbF}(0,T;\dbR^m)$, the MF-SDE \rf{fmf} with the initial condition $x(0)=x_0$ admits a unique strong solution $x(\cd)\in L^2_{\dbF}(\Omega;C([0,T];\dbR^n))$; see, e.g., \cite[Proposition 2.6]{Yong13}, \cite[Theorem 4.21]{Carmona-Delarue18i}.
To indicate the dependence on the initial state $x(0)=x_0$ and the control $u(\cd)$, we write the state as $x(\cd;x_0,u(\cd))$. If no confusion arises, we write it as $x(\cd)$ for short.

Throughout this paper, we adopt the following notions of controllability in the $L^2$ sense.

\begin{definition}\label{def-control}
The MF-SDE system \eqref{fmf} is said to be terminally controllable on $[0,T]$ if, for any $x_T\in L^2_{\cF_T}(\Omega;\dbR^n)$, there exists a pair $(x_0,u(\cd))\in \dbR^n\times L^2_\dbF(0,T;\dbR^m)$ such that the corresponding state satisfies $x(T;x_0, u(\cd))=x_T\;\as$
If \eqref{fmf} is terminally controllable on $[0,T]$ for every $T\in(0,\infty)$, then it is simply said to be terminally controllable.

The MF-SDE system \eqref{fmf} is said to be exactly controllable on $[0,T]$ if, for any initial state $x_0\in\dbR^n$ and any target $x_T\in L^2_{\cF_T}(\Omega;\dbR^n)$,  there exists a control $u(\cd)\in L^2_\dbF(0,T;\dbR^m)$ such that the corresponding state satisfies $x(T;x_0, u(\cd))=x_T \;\as$
If \eqref{fmf} is exactly controllable on $[0,T]$ for every $T\in(0,\infty)$, then it is simply said to be exactly controllable.
\end{definition}

The following notions concern optimal controls of NP problems.
\begin{definition}\label{def-minimal-norm}
For any $T>0$, $x_0\in\dbR^n$, and $x_T\in L^2_{\mf_T}(\Omega;\dbR^n)$, 
a control $u^*(\cdot)\in\cU(T,x_0,x_T)$ is called a minimal norm control of ${\rm\bf (NP)}^T_{x_0,x_T}$ if
$\cJ_F\big(u^*(\cdot)\big)
=\inf_{u(\cdot)\in\cU(T,x_0,x_T)}\cJ_F\big(u(\cd)\big)$,
where $\cU(T,x_0,x_T)$ and $\cJ_F(\cd)$ are defined in \rf{ad-control} and \rf{cost-F}, respectively.
The corresponding state $x^*(\cd)\in L^2_\dbF(\Omega;C([0,T];\dbR^n))$ is called an optimal state, and $\big(x^*(\cd),u^*(\cd)\big)$ is called an optimal pair.
\end{definition}

The following result will be used in Section \ref{sec-representation} to decompose the control variable under the full-row-rank condition for $D$ and $\widehat D=D+\bar D$.

\begin{lemma}\label{wang922l1}
Let $K\in\dbR^{n\times m}$ be of full row rank and set $r\deq m-n$. Define
\bel{wang922e2}
K^\dagger= K^\top(KK^\top)^{-1},
\qquad
\Pi_K= I_m-K^\dagger K.
\ee
Then $K^\dagger K$ is the orthogonal projection of $\dbR^m$ onto $\Ran(K^\top)$, while $\Pi_K$ is the orthogonal projection onto $\Ker K$.
Moreover, there exists a matrix $U_K\in\dbR^{m\times r}$ such that
\bel{wang922e3}
KU_K=0,
\qquad
U_K^\top U_K=I_r,
\qquad
U_KU_K^\top=\Pi_K.
\ee
\end{lemma}

\begin{proof}
Since $K$ has full row rank, $KK^\top$ is positive definite and hence invertible. Set $P_K= K^\dagger K=K^\top(KK^\top)^{-1}K$.
Then $P_K^\top=P_K$, and $P_K^2
=P_K$.
Hence $P_K$ is an orthogonal projection; see, e.g., \cite[Lecture 6]{Trefethen-Bau97}. 

Clearly, $\Ran(P_K)\subseteq\Ran(K^\top)$. 
Conversely, for any $y\in\Ran(K^\top)$, there exists $z\in\dbR^n$ such that $y=K^\top z$.
Hence, $
P_Ky
=K^\top(KK^\top)^{-1}KK^\top z
=K^\top z
=y$.
Therefore,
\[
\Ran(P_K)=\Ran(K^\top),
\]
and $P_K=K^\dagger K$ is the orthogonal projection onto $\Ran(K^\top)$.
On the other hand, since $\Ran(K^\top)=(\Ker K)^\perp$, it follows that
$\Pi_K=I_m-P_K$ is the orthogonal projection onto $\Ker K$.

Since $\Rank K=n$, it follows that $\dim(\Ker K)=m-n=r$.
Let $\{u_1,\cds,u_r\}$ be an orthonormal basis of $\Ker K$, and set
\[U_K=(u_1,\cds,u_r)\in\dbR^{m\times r}.\]
Then
\[
KU_K=0,
\qquad
U_K^\top U_K=I_r,
\qquad
\Ran(U_K)=\Ker K.
\]
Consequently,
$U_KU_K^\top$ is the orthogonal projection onto $\Ran(U_K)=\Ker K$. By the uniqueness of the orthogonal projection onto a given subspace,
\[
U_KU_K^\top=\Pi_K.
\]
This proves the assertion \rf{wang922e3} and completes the proof.
%
\end{proof}

\section{An equivalent rank condition for terminal controllability}
\label{sec-rank}

To study Problem {\bf (NP)$_{x_0,x_T}^T$} for arbitrary $T\in(0,\infty)$, $x_0\in\dbR^n$, and $x_T\in L^2_{\cF_T}(\Omega;\dbR^n)$, one needs to ensure that the admissible control set $\cU(T,x_0,x_T)$ is nonempty for any triple $(T,x_0,x_T)$. This naturally motivates the study of exact controllability for the underlying MF-SDE system. In this section, we establish an equivalent rank condition for the terminal controllability of \rf{fmf}; see Theorem \ref{rank-condition}. In particular, this condition is necessary for the exact controllability of \rf{fmf}.

Define $\cT: \mathbb R^n\times L^2_{\mathbb F}(0,T;\mathbb R^m) \longrightarrow L^2_{\mathcal F_T}(\Omega;\mathbb R^n)$ by $\cT (x_0,u(\cd))=x(T;x_0,u(\cd))$.
By the linearity and stability of the MF system \rf{fmf}, $\cT$ is a bounded linear operator. 
The following lemma provides a uniform estimate associated with the terminal controllability of \rf{fmf}, which will be used to prove the necessity of the rank condition in Theorem \ref{rank-condition}.

\begin{lemma}\label{Wang921l1}
Suppose that the MF-SDE system \rf{fmf} is terminally controllable on $[0,T]$.
Then there exists a constant $\cC>0$ such that, for every $x_T\in L^2_{\mathcal F_T}(\Omega;\mathbb R^n)$ satisfying
$\|x_T\|_{L^2_{\mathcal F_T}(\Omega;\mathbb R^n)}\leq 1$,
there is a pair
$(x_0,u(\cd))\in \mathbb R^n\times L^2_{\mathbb F}(0,T;\mathbb R^m)$
such that $x(T;x_0,u(\cd))=x_T$, and
\begin{align}\label{Wang921e1}
\|x_0\|^2+\me\int_0^T\|u(t)\|^2\rd t \leq \cC.
\end{align}
Moreover, the corresponding state satisfies
\begin{align}\label{Wang921e2}
\me\Big[\sup_{t\in[0,T]}\big\|x\big(t;x_0,u(\cd)\big)\big\|^2\Big]\leq \cC.
\end{align}
\end{lemma}

\begin{proof}
Set $X\deq\mathbb R^n\times L^2_{\mathbb F}(0,T;\mathbb R^m)$,
$Y\deq L^2_{\mathcal F_T}(\Omega;\mathbb R^n)$,
and equip $X$ with the norm
$\|(x_0,u(\cd))\|_X^2\deq \|x_0\|^2+\|u(\cd)\|^2_{L^2_\dbF(0,T;\dbR^m)}$.
Since \rf{fmf} is terminally controllable on $[0,T]$, $\cT$ is surjective. Hence, by the open mapping theorem, $\cT(B_X(0,1))$ is an open neighborhood of $0$ in $Y$.
Therefore, there exists a constant $r>0$ such that
$B_Y(0,r)\subseteq\cT(B_X(0,1))$.
Let $x_T\in Y$ satisfy $\|x_T\|_{Y}\leq1$.
Then
\[
\big\|r\frac {x_T} 2 \big\|_{Y}\leq\frac r2<r.
\]
Consequently, there exists $(z_0,v(\cd))\in B_X(0,1)$ such that
$\cT(z_0,v(\cd))=r\frac {x_T} 2$.
By the linearity of $\cT$, setting $(x_0,u(\cd))\deq\frac2r(z_0,v(\cd))$
gives
$\cT(x_0,u(\cd))=x(T;x_0,u(\cd))=x_T$.
Moreover,
\[
\|x_0\|^2+\|u(\cd)\|^2_{L^2_\dbF(0,T;\dbR^m)}
=\frac4{r^2}\|(z_0,v(\cd))\|_X^2
\leq\frac4{r^2}.
\]
Thus the assertion \rf{Wang921e1} follows by taking $\cC=4r^{-2}$.

Finally, by the stability estimate for the MF-SDE system \rf{fmf}, there exists a constant $\cC>0$, depending only on the coefficients of \rf{fmf} and $T$, such that
\[
\me\Big[\sup_{t\in[0,T]}\big\|x\big(t;x_0,u(\cd)\big)\big\|^2\Big]
\leq
\cC\Big[
\|x_0\|^2+\me\int_0^T\|u(t)\|^2\rd t
\Big].
\]
Combining this estimate with \rf{Wang921e1} yields \rf{Wang921e2}. This completes the proof.
\end{proof}

\begin{theorem}\label{rank-condition}
For any $T>0$, the MF-SDE system \rf{fmf} is terminally controllable on $[0,T]$ if and only if
\bel{rank-conditions}
\Rank(D)=\Rank(D+\bar D)=n.
\ee
Consequently, \rf{fmf} is terminally controllable if and only if \rf{rank-conditions} holds.
%
\end{theorem}

\begin{proof}
The sufficiency follows from \cite[Proposition 5]{Goreac14}. It therefore remains to prove the necessity. 
We establish the two rank conditions separately by contradiction.

\noindent
\textbf{Step 1: Full row rank of $D$.}
Suppose, to the contrary, that $\Rank D<n$. Then we can find a unit vector $q\in\dbR^n$ such that 
\bel{wang921e5}
q^\top D=0.
\ee 

For any $T\in (0,\infty)$ and any $\e\in(0,T)$, set 
\bel{wang921e3}
\th_\e=\frac{W(T-\e)}{\sqrt{T-\e}}
\qq\mbox{and}\qq
\xi_\e=\e^{-1/2}q\, \th_\e\big[W(T)-W(T-\e)\big].
\ee
Since $\th_\e$ is $\mf_{T-\e}$-measurable and $W(T)-W(T-\e)$ is independent of $\mf_{T-\e}$, we have
\bel{wang921e4}
\me\big[\th_\e\big]=0,\qq
\me\big[\|\xi_\e\|^2\big]=1,\qq
\me\big[\xi_\e\,\big|\,\mf_{T-\e}\big]=0.
\ee
Since \rf{fmf} is terminally controllable on $[0,T]$, Lemma \ref{Wang921l1}, applied to the target $\xi_\e\in L^2_{\mf_{T}}(\Omega;\dbR^n)$, ensures the existence of $(x_{0,\e}, u_\e(\cd))\in \dbR^n\times L^2_\dbF(0,T;\dbR^m)$ such that
\bel{wang921e8}
x_\e(T)=\xi_\e, \qq
\|x_{0,\e}\|^2+\me\int_0^T\|u_\e(t)\|^2\rd t\leq \cC,\qq
\me\Big[\sup_{t\in[0,T]}\|x_\e(t)\|^2\Big]\leq \cC,
\ee
where $x_\e(\cd)=x\big(\cd;x_{0,\e},u_\e(\cd)\big)$, and $\cC>0$ is independent of $\e$.

On the other hand, by \rf{wang921e5} and the decomposition $u_\e(\cd)=\wt u_\e(\cd)+\bar u_\e(\cd)$, we have $q^\top Du_\e(\cd)+q^\top \bar D \bar u_\e(\cd)=q^\top D\wt u_\e(\cd)+q^\top \h D \bar u_\e(\cd)=q^\top \h D \bar u_\e(\cd)$.
Define
\begin{align*}
&f_\e(\cd)=q^\top\big[Ax_\e(\cd)+\bar A\bar x_\e(\cd)
 +B u_\e(\cd)+\bar B\bar u_\e(\cd)\big],\\
\mbox{and}\q&
g_\e(\cd)=q^\top\big[Cx_\e(\cd)+\bar C\bar x_\e(\cd)\big].
\end{align*}
Then, integrating \rf{fmf} over $[T-\e,T]$ and using $x_\e(T)=\xi_\e$, we obtain
\begin{align}\label{wang921e6}
q^\top \xi_\e=q^\top x_\e(T-\e)+\int_{T-\e}^T f_\e(t)\rd t+\int_{T-\e}^T \big[g_\e(t)+q^\top \h D\bar u_\e(t)\big]\rd W(t).
\end{align}
Taking conditional expectation on both sides with respect to $\mf_{T-\e}$ and using \rf{wang921e4} gives
\begin{align}\label{wang921e7}
0=q^\top x_\e(T-\e)+\me\Big[\int_{T-\e}^T f_\e(t)\rd t \,\Big|\,\mf_{T-\e}\Big].
\end{align}
Subtracting \rf{wang921e7} from \rf{wang921e6} yields
\begin{align}
&q^\top \xi_\e-\int_{T-\e}^T q^\top \h D\bar u_\e(t)\rd W(t)\notag\\
&=\Big\{\int_{T-\e}^T f_\e(t)\rd t-\me\Big[\int_{T-\e}^T f_\e(t)\rd t \,\Big|\,\mf_{T-\e}\Big] \Big\}+\int_{T-\e}^T g_\e(t)\rd W(t)\label{wang921e10}\\
&=: I_1(\e)+I_2(\e).\notag
\end{align}

We next estimate $I_1(\e)$ and $I_2(\e)$ on the right-hand side of \rf{wang921e10}. For $I_1(\e)$, the $L^2$-contraction property of conditional expectation, H\"older's inequality, and \rf{wang921e8} imply that
\begin{align}
\me\big[\|I_1(\e)\|^2\big]
&\leq \e\me\int_{T-\e}^T\|f_\e(t)\|^2\rd t\notag\\
&\leq \cC\e\me\int_{T-\e}^T\big[\|x_\e(t)\|^2+\|u_\e(t)\|^2\big]\rd t\label{wang921e9}\\
&\leq \cC\e.\notag
\end{align}
For $I_2(\e)$, the It\^o isometry and \rf{wang921e8} lead to
\begin{align}\label{wang921e12}
\me\big[\|I_2(\e)\|^2\big]
\leq \cC\me\int_{T-\e}^T\sup_{s\in[0,T]}\|x_\e(s)\|^2\rd t
\leq \cC\e.
\end{align} 

It remains to obtain a lower bound for the left-hand side of
\rf{wang921e10}.
Since $\th_\e$ is $\mf_{T-\e}$-measurable, the conditional It\^o isometry gives
\begin{align*}
&\me\Big[q^\top \xi_\e \int_{T-\e}^T q^\top \h D\bar u_\e(t)\rd W(t)\Big]\\
&=\me\Big[\e^{-1/2} \th_\e \me\Big( \int_{T-\e}^T 1 \rd W(t) \int_{T-\e}^T q^\top \h D\bar u_\e(t)\rd W(t)\,\Big|\,\mf_{T-\e}\Big)\Big]\\
&=\me\Big[\e^{-1/2} \th_\e  \int_{T-\e}^T  q^\top \h D\bar u_\e(t)\rd t\Big].
\end{align*}
Since $\bar u_\e(\cd)$ is deterministic, it follows that
\[
\me\Big[q^\top \xi_\e \int_{T-\e}^T q^\top \h D\bar u_\e(t)\rd W(t)\Big]
=\e^{-1/2}\me\big[\th_\e\big] \int_{T-\e}^T  q^\top \h D\bar u_\e(t)\rd t
=0,
\]
where the last equality follows from \rf{wang921e4}.
Thus,
\begin{align}\label{wang921e13}
\me\Big[\Big\|q^\top \xi_\e\!-\! \int_{T-\e}^T \!q^\top \h D\bar u_\e(t)\rd W(t)\Big\|^2\Big]
=\me\big[\|\xi_\e\|^2\big]+\int_{T-\e}^T\!\big\| q^\top \h D\bar u_\e(t) \big\|^2\rd t
\geq 1.
\end{align}

Finally, combining \rf{wang921e10}, \rf{wang921e9}, \rf{wang921e12}, and \rf{wang921e13}, we conclude that
\beq
1\leq 2\me\big[\|I_1(\e)\|^2+\|I_2(\e)\|^2\big]\leq \cC\e.
\eeq
Since $\e>0$ is arbitrary, letting $\e\downarrow0$ yields a contradiction. Therefore, $\Rank(D)=n$.

\noindent
\textbf{Step 2: Full row rank of $\h D=D+\bar D$.}
Suppose that $\Rank(\h D)<n$, and choose a unit vector $q_0 \in\dbR^n$ satisfying $q_0^\top \h D=0$. 
For any $\e\in(0,T)$, set the target $\eta_\e=\e^{-1/2}q_0\big[W(T)-W(T-\e)\big]$. It follows that $\me\big[\|\eta_\e\|^2\big]=1$ and
$\me\big[\eta_\e\,\big|\,\mf_{T-\e}\big]=0$. 
Terminal controllability of \rf{fmf} on $[0,T]$, together with Lemma \ref{Wang921l1}, provides $(z_{0,\e}, v_\e(\cd))\in \dbR^n\times L^2_\dbF(0,T;\dbR^m)$ satisfying 
$X_\e(T)=\eta_\e$, 
$\|z_{0,\e}\|^2+\me\int_0^T\|v_\e(t)\|^2\rd t\leq \cC$, and
$\me\big[\sup_{t\in[0,T]}\|X_\e(t)\|^2\big]\leq \cC$,
where $X_\e(\cd)=x\big(\cd;z_{0,\e},v_\e(\cd)\big)$.

Since
$v_\e(\cd)=\wt v_\e(\cd)+\bar v_\e(\cd)$,
$\bar v_\e(\cd)=\me[v_\e(\cd)]$, and $q_0^\top\h D=0$, we have
$q_0^\top\big[Dv_\e(\cd)+\bar D\bar v_\e(\cd)\big]
=q_0^\top D\wt v_\e(\cd)+q_0^\top\h D\bar v_\e(\cd)
=q_0^\top D\wt v_\e(\cd)$.
Proceeding as in Step 1 to derive \rf{wang921e10}, we have
\begin{align*}
&q_0^\top \eta_\e-\int_{T-\e}^T q_0^\top D\wt v_\e(t)\rd W(t)\\
&\qq=\int_{T-\e}^T F_\e (t)\rd t-\me\Big[\int_{T-\e}^T F_\e(t)\rd t \,\Big|\,\mf_{T-\e}\Big]+\int_{T-\e}^T G_\e(t)\rd W(t),
\end{align*}
where $F_\e(\cd)=q_0^\top\big[AX_\e(\cd)+\bar A\bar X_\e(\cd)
 +B v_\e(\cd)+\bar B\bar v_\e(\cd)\big]$, and
$G_\e(\cd)=q_0^\top\big[CX_\e(\cd)+\bar C\bar X_\e(\cd)\big]$.
As in \rf{wang921e9} and \rf{wang921e12}, we have
\beq
\me\Big[\Big\|\int_{T-\e}^T F_\e (t)\rd t-\me\Big[\int_{T-\e}^T F_\e(t)\rd t \,\Big|\,\mf_{T-\e}\Big]+\int_{T-\e}^T G_\e(t)\rd W(t)\Big\|^2\Big]
\leq \cC\e.
\eeq
On the other hand,
\begin{align*}
&\me\Big[\Big\|q_0^\top \eta_\e-\int_{T-\e}^T q_0^\top D\wt v_\e(t)\rd W(t)\Big\|^2\Big]\\
&=1+\me\int_{T-\e}^T\big\| q_0^\top D\wt v_\e(t)\big\|^2\rd t
-2 \me\int_{T-\e}^T\e^{-1/2}q_0^\top D \me\big[\wt v_\e(t)\big]\rd t\\
&\geq 1.
\end{align*}
Here we apply $\me\big[\wt v_\e(t)\big]=\me\big[v_\e(t)-\me\big[v_\e(t)]\big]=0$.
Combining the above estimates yields $1\leq \cC\e$ for any $\e\in (0,T)$, which is impossible as $\e\downarrow0$. Hence, $\Rank(D+\bar D)=n$.
This completes the proof.
\end{proof}

\section{Representation of the minimal norm control}\label{sec-representation}

In this section, we derive a Riccati-based representation of the minimal norm control for Problem {\bf (NP)$_{x_0,x_T}^T$} under the exact controllability assumption on the underlying MF-SDE system \rf{fmf}.
For simplicity of notation, in what follows we shall often suppress the time variable when no confusion can arise.

\subsection{Equivalent transformation of the minimal norm problem}\label{subsec-transformation}

We first transform the original MF-SDE system \rf{fmf} into an equivalent backward system. Suppose that
\[
\Rank D=\Rank\h D=n,
\qquad
\h D=D+\bar D.
\]
By Theorem \ref{rank-condition}, this rank condition is automatically satisfied whenever \rf{fmf} is exactly controllable on $[0,T]$. In particular, $DD^\top$ and $\h D\h D^\top$ are both invertible.


Applying Lemma \ref{wang922l1} with $K=D$ and $K=\h D$, respectively, and setting $r\deq m-n$, we fix matrices $U,\h U\in\dbR^{m\times r}$ whose columns form orthonormal bases of $\Ker D$ and $\Ker\h D$, respectively. Accordingly, 
\begin{align}
&D^\dagger\deq D^\top(DD^\top)^{-1},
\qquad
\Pi_D\deq I_m-D^\dagger D,\notag\\
&\h D^\dagger\deq \h D^\top(\h D\h D^\top)^{-1},
\qquad
\Pi_{\h D}\deq I_m-\h D^\dagger\h D,\notag\\
&DU=0,\qquad U^\top U=I_r,\qquad UU^\top=\Pi_D,\notag\\
&\h D\h U=0,\qquad \h U^\top\h U=I_r,\qquad
\h U\h U^\top=\Pi_{\h D}.
\label{decom}
\end{align}
We now introduce the transformed variables
\bel{wang922e9}
\left\{
\bal
&z(\cd)=C\wt x(\cd)+\h C\bar x(\cd)
+D\wt u(\cd)+\h D\bar u(\cd),\\
&v(\cd)=U^\top\wt u(\cd)+\h U^\top\bar u(\cd).
\eal
\right.
\ee
Since $\me[\wt u(\cd)]=0$, it follows from \rf{wang922e9} that
\begin{align*}
&\wt z(\cd)-C\wt x(\cd)=D\wt u(\cd),
&&\bar z(\cd)-\h C\bar x(\cd)=\h D\bar u(\cd),\\
\mbox{and }\qq&
\wt v(\cd)=U^\top\wt u(\cd),
&&\bar v(\cd)=\h U^\top\bar u(\cd).
\end{align*}
Hence, by \rf{decom}, the control $u(\cd)$ can be represented in terms of $x(\cd)$, $z(\cd)$, and $v(\cd)$ as
\bel{wang922e1}
\bal
u(\cd)&=D^\dagger D\wt u(\cd)+\Pi_D \wt u(\cd) + \h D^\dagger \h D\bar u(\cd)+\Pi_{\h D} \bar u(\cd)\\
&=D^\dagger \big[\wt z(\cd)-C\wt x(\cd)\big]+U \wt v(\cd) 
+ \h D^\dagger \big[\bar z(\cd)-\h C \bar x(\cd)\big]+\h U \bar v(\cd)\\
&=: f\big(x(\cd),z(\cd),v(\cd)\big).\\
\eal
\ee
Set
\bel{wang922e5}
\bal
&\BA=A-BD^\dagger C,\q
\BK=BD^\dagger,\q
\BL=BU,\\
\mbox{and}\q&
\h\BA=\h A-\h B\h D^\dagger \h C,\q
\h \BK=\h B\h D^\dagger,\q
\h \BL=\h B\h U.
\eal
\ee
Substituting \rf{wang922e1} into \rf{fmf}, we obtain the equivalent backward system
\bel{fmf-2}
\bal
& \mathrm dx(t)=\big[\BA \wt x+\h \BA \bar x
+\BK \wt z+\h \BK\bar z+ \BL \wt v+\h \BL\bar v\big]\rd t
+z\rd W(t) \q t\in[0,T]. \\
\eal
\ee
In what follows, \rf{fmf-2} is considered together with the terminal condition $x(T)=x_T$ and is therefore referred to as a controlled MF-BSDE system, with $(x(\cd),z(\cd))$ as the solution pair and $v(\cd)\in L^2_\dbF(0,T;\dbR^r)$ as the control, where $z(\cd)$ is the martingale integrand.
To emphasize the dependence on the terminal state $x_T$ and the control $v(\cd)$, we write the first component of the corresponding solution as $x(\cd;x_T,v(\cd))$, or simply as $x(\cd)$ when no confusion can arise.

For any $x_0\in\dbR^n$ and $x_T\in L^2_{\mf_T}(\Omega;\dbR^n)$, define the admissible control set for the MF-BSDE system \rf{fmf-2} by
\bel{wang922e7}
\cV(T,x_0,x_T)
\deq
\big\{
v(\cdot)\in L^2_{\dbF}(0,T;\dbR^r)
\,\big|\,
x\big(0;x_T,v(\cdot)\big)=x_0
\big\}.
\ee
We say that the MF-BSDE system \rf{fmf-2} is exactly controllable on $[0,T]$ if the admissible control set $\cV(T,x_0,x_T)\neq\emptyset$ for every $x_0\in\dbR^n$ and every $x_T\in L^2_{\mf_T}(\Omega;\dbR^n)$.
If \rf{fmf-2} is exactly controllable on $[0,T]$ for every $T\in(0,\infty)$, then it is simply said to be exactly controllable.

We next introduce a constrained optimal control problem associated with \rf{fmf-2}. Define
\bel{wang922e10}
\bal
&Q=C^\top (DD^\top)^{-1} C,
\qquad
S=-(DD^\top)^{-1} C,
\qquad
R=(DD^\top)^{-1},\\
\mbox{and}\q&\widehat Q=\widehat C^\top(\widehat D\widehat D^\top)^{-1}\widehat C,
\qquad
\widehat S=-(\widehat D\widehat D^\top)^{-1}\widehat C,
\qquad
\widehat R=(\widehat D\widehat D^\top)^{-1}.
\eal
\ee
Using \rf{decom}, in particular, $D U=\h D \h U=O_{n\times r}$, and $U^\top U=\h U^\top \h U=I_r$, we arrive at
\bel{cost-functional}
\bal
\cJ_B\big(v(\cd)\big)
&\deq\big\|f\big(x(\cd),z(\cd),v(\cd)\big)\big\|^2_{L^2_\dbF(0,T;\dbR^m)}\\
&=\me \int_0^T\Big[\lan Q \wt x,\wt x\ran+2\lan S\wt x, \wt z\ran+\lan R\wt z, \wt z\ran +\|\wt v\|^2
\Big]\rd t\\
&\q+\int_0^T\Big[\lan \h Q \bar x,\bar x\ran+2\lan \h S\bar x, \bar z\ran+\lan \h R\bar z, \bar z\ran +\|\bar v\|^2
\Big]\rd t.
\eal
\ee
Then, the constrained optimal control problem associated with \rf{fmf-2} is formulated as follows.\\
{\bf Problem (COP)$_{x_0,x_T}^T$}. For any $x_0\in \dbR^n$ and any $x_T\in L^2_{\mf_T}(\Omega;\dbR^n)$, find a control $v^*(\cd)\in \cV(T,x_0,x_T)$ minimizing the cost functional $\cJ_B\big(v(\cd)\big)$.

The following proposition establishes the equivalence between the exact controllability of \rf{fmf} and that of \rf{fmf-2}, as well as the equivalence between Problems {\bf (NP)$_{x_0,x_T}^T$} and {\bf (COP)$_{x_0,x_T}^T$}.

\bp{wang922p1}
Suppose that $\Rank(D)=\Rank(D+\bar D)=n$. Then the exact controllability of the MF-SDE system \rf{fmf} is equivalent to that of the MF-BSDE system \rf{fmf-2}.
Furthermore, suppose that \rf{fmf} is exactly controllable on $[0,T]$. Then, for every $x_0\in\dbR^n$ and $x_T\in L^2_{\mf_T}(\Omega;\dbR^n)$, Problems {\rm\bf (NP)$_{x_0,x_T}^T$} and {\rm\bf (COP)$_{x_0,x_T}^T$} admit unique optimal controls $u^*(\cd)$ and $v^*(\cd)$, respectively, and
\bel{wang922e8}
u^*(\cd)= f\big(x^*(\cd),z^*(\cd),v^*(\cd)\big),
\ee
where $\big(x^*(\cd),z^*(\cd)\big)$ is the solution to the MF-BSDE \rf{fmf-2} with terminal condition $x^*(T)=x_T$ and control $v^*(\cd)$.
\ep

\begin{proof}
We first prove the equivalence of exact controllability. Let $u(\cd)$ be any control for \rf{fmf}, with $u(\cd)\in L^2_\dbF(0,T;\dbR^m)$, and let $x(\cd)$ be the corresponding state. Define $z(\cd)$ and $v(\cd)$ by \rf{wang922e9}. Then \rf{wang922e1} and \rf{wang922e5} imply that $\big(x(\cd),z(\cd)\big)$ satisfies \rf{fmf-2}. Conversely, for any $v(\cd)\in L^2_\dbF(0,T;\dbR^r)$ and terminal state $x_T$, let $\big(x(\cd),z(\cd)\big)$ be the solution of \rf{fmf-2}, and define $u(\cd)$ by \rf{wang922e1}. A direct substitution shows that $x(\cd)$ satisfies \rf{fmf}. Hence, for every $T>0$, the two systems have the same reachable initial-terminal pairs, and their exact controllability on $[0,T]$ is equivalent. This argument is analogous to that in \cite[Theorem 3.6]{Yu21}.

Now suppose that \rf{fmf} is exactly controllable on $[0,T]$. Define the bounded linear terminal map $\cT_0:L^2_\dbF(0,T;\dbR^m)\longrightarrow L^2_{\mf_T}(\Omega;\dbR^n)$ by $\cT_0u(\cd)= x(T;0,u(\cd))$.
For fixed $x_0\in\dbR^n$, the linearity of \rf{fmf} gives $x(T;x_0,u(\cd))=x(T;x_0,0)+\cT_0u(\cd)$.
Hence,
\[
\cU(T,x_0,x_T)
=\cT_0^{-1}\big(\{x_T-x(T;x_0,0)\}\big),
\]
which is closed in $L^2_\dbF(0,T;\dbR^m)$. Moreover, it is an affine subset, and the exact controllability of \rf{fmf} guarantees that it is nonempty. Therefore, $\cU(T,x_0,x_T)$ is a nonempty closed affine subset of the Hilbert space $L^2_\dbF(0,T;\dbR^m)$.
Hence, by the projection theorem in Hilbert spaces, there exists a unique $u^*(\cd)\in\cU(T,x_0,x_T)$ such that $\cJ_F\big(u^*(\cd)\big)=\inf_{u(\cd)\in \cU(T,x_0,x_T)}\cJ_F\big(u(\cd)\big)$. Therefore, Problem {\bf (NP)$_{x_0,x_T}^T$} admits a unique optimal control $u^*(\cd)$.

By the first assertion, \rf{fmf-2} is also exactly controllable on $[0,T]$, and therefore $\cV(T,x_0,x_T)$ is a nonempty closed affine subset of $L^2_\dbF(0,T;\dbR^r)$. Moreover, the matrices
\begin{align}\label{wang923e1}
\Theta\deq
\begin{pmatrix}
Q&S^\top\\
S&R
\end{pmatrix},
\qquad
\h\Theta\deq
\begin{pmatrix}
\h Q&\h S^\top\\
\h S&\h R
\end{pmatrix}
\end{align}
are positive semidefinite. Indeed,
$\Theta=\big(C\,-I_n\big)^\top(DD^\top)^{-1}\big(C\,-I_n\big)\geq0$,
and similarly $\h\Theta\geq0$. Hence, by \rf{cost-functional},
\begin{align*}
\cJ_B\big(v(\cd)\big)
\geq \me\int_0^T\|\wt v(t)\|^2\rd t +\int_0^T\|\bar v(t)\|^2\rd t
=\|v(\cd)\|^2_{L^2_\dbF(0,T;\dbR^r)}.
\end{align*}
Therefore, $\cJ_B(\cd)$ is coercive. Since the solution $\big(x(\cd),z(\cd)\big)$ of \rf{fmf-2} depends affinely and continuously on $v(\cd)$, the functional $\cJ_B(\cd)$ is continuous and strictly convex on $L^2_\dbF(0,T;\dbR^r)$. In particular, it is weakly lower semicontinuous. Therefore, by the direct method of the calculus of variations, Problem {\bf (COP)$_{x_0,x_T}^T$} admits a unique optimal control $v^*(\cd)$.

It remains to establish \rf{wang922e8}. For any $u(\cd)\in\cU(T,x_0,x_T)$, define $z(\cd)$ and $v(\cd)$ by \rf{wang922e9}. Then, by \rf{wang922e1}, the corresponding $(x(\cd),z(\cd))$ satisfies \rf{fmf-2}, $v(\cd)\in\cV(T,x_0,x_T)$, and $u(\cd)=f\big(x(\cd),z(\cd),v(\cd)\big)$.
Conversely, for any control $v(\cd)\in\cV(T,x_0,x_T)$, if $(x(\cd),z(\cd))$ is the corresponding solution of \rf{fmf-2}, then \rf{wang922e1} defines a control $u(\cd)=f\big(x(\cd),z(\cd),v(\cd)\big)\in\cU(T,x_0,x_T)$.
Moreover, by \rf{cost-functional},
\[
\cJ_F\big(u(\cd)\big)=\cJ_B\big(v(\cd)\big).
\]
Thus, Problems {\bf (NP)$_{x_0,x_T}^T$} and {\bf (COP)$_{x_0,x_T}^T$} are equivalent under the transformations \rf{wang922e9} and \rf{wang922e1}. 
Let $v^*(\cd)$ be the optimal control of {\bf (COP)$_{x_0,x_T}^T$}, and let $\big(x^*(\cd),z^*(\cd)\big)$ be the corresponding solution of \rf{fmf-2}. Set $\widehat u(\cd)\deq f\big(x^*(\cd),z^*(\cd),v^*(\cd)\big)$.
Then $\widehat u(\cd)\in\cU(T,x_0,x_T)$. For any $u(\cd)\in\cU(T,x_0,x_T)$, let $v(\cd)\in\cV(T,x_0,x_T)$ be the control defined by \rf{wang922e9}. By \rf{cost-functional} and the optimality of $v^*(\cd)$,
\[
\cJ_F\big(\widehat u(\cd)\big)
=\cJ_B\big(v^*(\cd)\big)
\leq\cJ_B\big(v(\cd)\big)
=\cJ_F\big(u(\cd)\big).
\]
Hence $\widehat u(\cd)$ is an optimal control of {\bf (NP)$_{x_0,x_T}^T$}. By the uniqueness of the latter,
\beq
u^*(\cd)=f\big(x^*(\cd),z^*(\cd),v^*(\cd)\big),
\eeq
which proves \rf{wang922e8}.
This completes the proof.
\end{proof}

As an immediate consequence of Theorem \ref{rank-condition} and Proposition \ref{wang922p1}, we obtain the following characterization of exact controllability in terms of the solvability of the minimal norm control problem.

\bc{control=NP}
The MF-SDE system \rf{fmf} is exactly controllable on $[0,T]$ if and only if, for every $x_0\in\dbR^n$ and $x_T\in L^2_{\mf_T}(\Omega;\dbR^n)$, Problem {\rm\bf (NP)$_{x_0,x_T}^T$} admits a unique minimal norm control.
\ec

\begin{proof}
Suppose first that \rf{fmf} is exactly controllable on $[0,T]$. Then Proposition \ref{wang922p1} implies that, for every $x_0\in\dbR^n$ and $x_T\in L^2_{\mf_T}(\Omega;\dbR^n)$, Problem {\rm\bf (NP)$_{x_0,x_T}^T$} admits a unique optimal control, which is precisely its unique minimal norm control.

Conversely, suppose that, for every $x_0\in\dbR^n$ and $x_T\in L^2_{\mf_T}(\Omega;\dbR^n)$, Problem {\rm\bf (NP)$_{x_0,x_T}^T$} admits a unique minimal norm control. Then $\cU(T,x_0,x_T)\neq\emptyset$ for every such pair. By Definition \ref{def-control}, \rf{fmf} is exactly controllable on $[0,T]$.
\end{proof}

\subsection{Parameterized optimal control problem}\label{subsec-para}

To handle the initial constraint $x(0)=x_0$ in Problem {\bf (COP)$_{x_0,x_T}^T$}, we introduce a parameter $\l\in\dbR^n$ as a Lagrange multiplier and consider an associated unconstrained LQ problem. For any $\l,x_0\in\dbR^n$ and $x_T\in L^2_{\mf_T}(\Omega;\dbR^n)$, we formulate the following parameterized optimal control problem.\\
{\bf Problem (OP)$_{x_0,x_T,\l}^T$.} Search for a control $v_\l^*(\cd)\in L^2_\dbF(0,T;\dbR^r)$ minimizing the following $\l$-dependent cost functional
\[
\cJ_\l\big(v(\cd)\big)=\cJ_B\big(v(\cd)\big)+2\big\lan \l, x\big(0;x_T,v(\cd)\big)-x_0\big\ran,
\]
where $x\big(\cd;x_T,v(\cd)\big)$ is the state of the MF-BSDE \rf{fmf-2} with the terminal condition $x\big(T;x_T,v(\cd)\big)=x_T$.

The following result establishes the connection between Problems {\bf (OP)$_{x_0,x_T,\l}^T$} and {\bf (COP)$_{x_0,x_T}^T$}.

\bl{wang923l1}
For any $\l, x_0\in\dbR^n$ and $x_T\in L^2_{\mf_T}(\Omega;\dbR^n)$, Problem {\bf (OP)$_{x_0,x_T,\l}^T$} admits a unique optimal triple $\big(x^*_\l(\cd),z^*_\l(\cd),v^*_\l(\cd)\big)$. Moreover, there exists a co-state process $y^*_\l(\cd)$ such that
$\big(x^*_\l(\cd),y^*_\l(\cd),z^*_\l(\cd),v^*_\l(\cd)\big)$ satisfies
\bel{fbsde}
\left\{
\bal
&\mathrm{d}x^*_\l(t)\!=\!\big[\BA \wt x^*_\l\!+\!\h \BA \bar x^*_\l
\!+\!\BK \wt z^*_\l\!+\!\h \BK\bar z^*_\l\!+\! \BL \wt v^*_\l\!+\!\h \BL\bar v^*_\l\big]\rd t\!+\!z^*_\l\rd W(t),\q t\in [0,T],\\
&\mathrm{d}y^*_\l(t)\!=\!\big[\!-\!\BA^{\!\top} \wt y^*_\l\!-\!\h\BA^{\!\top}\bar y^*_\l\!+\!Q\wt x^*_\l\!+\!\h Q\bar x^*_\l\!+\!S^\top \wt z^*_\l\!+\!\h S^\top \bar z^*_\l\big]\rd t\\
&\qq\qq +\big[\!-\!\BK^{\!\top} \wt y^*_\l\!-\!\h\BK^{\!\top}\bar y^*_\l\!+\!S\wt x^*_\l\!+\!\h S\bar x^*_\l\!+\!R \wt z^*_\l\!+\!\h R\bar z^*_\l\big]\rd W(t), \q t\in [0,T],\\
&x^*_\l(T)\!=\!x_T,\q y^*_\l(0)\!=\!\l,
\eal
\right.
\ee
together with the optimality condition
\bel{optimality}
\wt v^*_\l=\BL^\top \wt y^*_\l,\q \bar v^*_\l=\h\BL^\top \bar y^*_\l \q \ae
\ee
Furthermore, if $x^*_\l(0)=x_0$ for some $\l\in\dbR^n$, then $v^*_\l(\cd)$ is the unique optimal control of Problem {\bf (COP)$_{x_0,x_T}^T$}.
\el

\begin{proof}
Since $\Theta$ and $\h\Theta$ defined in \rf{wang923e1} are positive semidefinite and the control weight in \rf{cost-functional} is the identity matrix, the quadratic part of $\cJ_\l$ is uniformly convex with respect to $v(\cd)$. 
The additional term $2\big\lan\l,x\big(0;x_T,v(\cd)\big)-x_0\big\ran$ is affine and continuous in $v(\cd)$. Hence Problem {\bf (OP)$_{x_0,x_T,\l}^T$} admits a unique optimal control $v^*_\l(\cd)$; the associated solution pair $\big(x^*_\l(\cd),z^*_\l(\cd)\big)$ is uniquely determined by \rf{fmf-2}. 
The optimality system \rf{fbsde}--\rf{optimality} follows from the stochastic maximum principle; see, e.g., \cite[Theorem 3.1]{Li-Sun-Xiong19}.

Now suppose that $x^*_\l(0)=x_0$ for some $\l\in\dbR^n$. Then $v^*_\l(\cd)\in\cV(T,x_0,x_T)$.
For any control $v(\cd)\in \cV(T,x_0,x_T)\subseteq L^2_\dbF(0,T;\dbR^r)$, it follows that $x\big(0;x_T,v(\cd)\big)=x_0$. Hence, by the definition of $\cJ_\l$ and the condition $x^*_\l(0)=x_0$, 
\[
\cJ_B\big(v^*_\l(\cd)\big)=\cJ_\l\big(v^*_\l(\cd)\big)
\leq \cJ_\l\big(v(\cd)\big)=\cJ_B\big(v(\cd)\big).
\]
Since $v(\cd)\in\cV(T,x_0,x_T)$ is arbitrary, it follows that $v^*_\l(\cd)$ is an optimal control of Problem {\bf (COP)$_{x_0,x_T}^T$}. Its uniqueness follows from the unique solvability of {\bf (COP)$_{x_0,x_T}^T$} established in Proposition \ref{wang922p1}.
\end{proof}

To derive an explicit representation of the optimal triple, we apply the Riccati decoupling technique to the optimality system \rf{fbsde}--\rf{optimality}. Motivated by the standard decoupling ansatz for LQ problems (see, e.g., \cite[Section 6.6]{Yong-Zhou99} and \cite[Section 3]{Li-Sun-Xiong19}), we seek
\begin{align}\label{feedback-x}
x^*_\l(\cd)=-P(\cd)\bar y^*_\l(\cd)-\Sigma(\cd) \wt y^*_\l(\cd)+\f(\cd).
\end{align}
This leads to the following two differential Riccati equations:
\begin{equation}\label{Riccati-Sigma}
\left\{
\begin{aligned}
&\dot \Sigma 
-\BA\Sigma -\Sigma \BA^\top-\Sigma Q\Sigma + \BL\BL^\top  \\
&\qq\q+( \BK+\Sigma  S^\top)(I+\Sigma R)^{-1}\Sigma ( \BK^\top+S\Sigma )=0,\q \mbox{on } [0,T],\\
&\Sigma (T)=0,
\end{aligned}
\right.
\end{equation}
and
\begin{equation}\label{Riccati-P}
\left\{
\begin{aligned}
&\dot P
-\h{\BA}P-P\h{\BA}^{\top}-P\h QP
+\h \BL\h \BL^{\top}\\
&\qq\q+(\h{\BK}+P\h S^\top)(I+\Sigma \h R)^{-1}
\Sigma (\h{\BK}^\top+\h S P)=0,\q\mbox{on }[0,T],\\
&P(T)=0.
\end{aligned}
\right.
\end{equation}
The remaining terms give a $(\Sigma(\cd),P(\cd))$-dependent MF-BSDE
\begin{equation}\label{bsde-varphi}
\left\{
\begin{aligned}
&\mathrm{d}\varphi(t)=
\Big[\cM_1\wt \varphi+ \h\cM_1\bar \varphi
+( \BK+\Sigma S^\top)(I+\Sigma R)^{-1}\wt\beta
\\
&\qq\qq+(\h{\BK}+P\h S^\top)(I+\Sigma \h R)^{-1}\bar\beta
\Big]\rd t
+\beta \rd W(t),\qq t\in[0,T],\\
&\varphi(T)=x_T,
\end{aligned}
\right.
\end{equation}
where
\bel{wang923e3}
\bal
&\cM_1=\BA+\Sigma Q-( \BK+\Sigma S^\top)(I+\Sigma R)^{-1}\Sigma S,\\
\mbox{and}\q &
\h \cM_1=\h{\BA}+P\h Q
-(\h{\BK}+P\h S^\top)(I+\Sigma \h R)^{-1}\Sigma \h S.
\eal
\ee

The following lemma establishes the well-posedness of \rf{Riccati-Sigma}--\rf{bsde-varphi} and provides the corresponding decoupled representations.

\bl{wang923l2}
The Riccati equations \rf{Riccati-Sigma} and \rf{Riccati-P} have unique solutions $\Sigma(\cd)\in C([0,T];\dbS^n_+)$ and $P(\cd)\in C([0,T];\dbS^n_+)$, respectively. Then the 
MF-BSDE \rf{bsde-varphi} admits a unique solution $\big(\f(\cd),\beta(\cd)\big)\in L^2_\dbF(\Omega;C([0,T];\dbR^n))\times L^2_\dbF(0,T;\dbR^n)$.
Moreover, the martingale integrand $z^*_\l(\cd)$ in \rf{fbsde} has the following representation:
\bel{feedback-z}
\bal
z^*_\l&=(I+\Sigma R)^{-1}\big[\Sigma  \big(\BK^\top+S\Sigma\big) \wt y^*_\l-\Sigma S\wt \varphi+\wt\beta\big]\\
&\q+(I+\Sigma \h R)^{-1}\big[\Sigma  \big(\h \BK^\top+\h S P\big) \bar y^*_\l-\Sigma \h S\bar \varphi+\bar\beta\big].
\eal
\ee
Furthermore, the co-state $y^*_\l(\cd)$ in \rf{fbsde} satisfies the following decoupled MF-SDE:
\begin{equation}\label{sde-y}
\left\{
\begin{aligned}
&\mathrm{d}y^*_\l(t)\!=\!
-\big[\cM_1^\top \wt y^*_\l\!+\!\h{\cM}_1^\top\bar y^*_\l
\!+\!\cM_2\wt\varphi\!+\!\h{\cM}_2\bar \varphi
\!+\!\cM_3\wt\beta\!+\!\h{\cM}_3\bar\beta
\big]\rd t\\
&\qq\qq
-\big[\cN_1 \wt y^*_\l\!+\!\h{\cN}_1\bar y^*_\l
\!+\!\cN_2\wt\varphi\!+\!\h{\cN}_2\bar \varphi
\!+\!\cN_3\wt\beta\!+\!\h{\cN}_3\bar\beta
\big]\rd W(t),\q t\in [0,T],\\
&y^*_\l(0)\!=\!\lambda,
\end{aligned}
\right.
\end{equation}
where
\bel{wang923e4}
\bal
&\cM_2
=-Q+S^\top(I+\Sigma R)^{-1}\Sigma S,\q 
\h{\cM}_2
=-\h Q+\h S^\top(I+\Sigma \h R)^{-1}\Sigma \h S,\\
&\cM_3
=-S^\top(I+\Sigma R)^{-1},\qquad\qq
\h{\cM}_3=-\h S^\top(I+\Sigma \h R)^{-1},\\
&\cN_1
=\BK^\top+S\Sigma  -R(I+\Sigma R)^{-1}\Sigma ( \BK^\top+S\Sigma ),\\
&\h{\cN}_1
=\h{\BK}^\top+\h SP
-\h R(I+\Sigma \h R)^{-1}\Sigma 
(\h{\BK}^\top+\h SP),\\
&\cN_2
=-S+R(I+\Sigma R)^{-1}\Sigma S,\qq
\h{\cN}_2
=-\h S+\h R(I+\Sigma \h R)^{-1}\Sigma \h S,\\
&\cN_3
=-R(I+\Sigma R)^{-1},\qquad\qq\q
\h{\cN}_3=-\h R(I+\Sigma \h R)^{-1}.
\eal
\ee

\el

\begin{proof}
The unique solvability of \rf{Riccati-Sigma} and \rf{Riccati-P} follows by an argument similar to that in \cite[Theorem 4.1]{Li-Sun-Xiong19}. Once $\Sigma(\cd)$ and $P(\cd)$ are determined, \rf{bsde-varphi} is a linear MF-BSDE with bounded deterministic coefficients. Hence, by \cite[Theorem 3.1]{Buckdahn-Li-Peng09}, it admits a unique solution.

We next derive \rf{feedback-z}. Applying It\^o's formula to \rf{feedback-x}  yields
\begin{align*}
\mathrm{d}x^*_\l(t)&=\cds\rd t+\Big\{\Sigma \big[\BK^\top \wt y^*_\l+\h\BK^\top\bar y^*_\l+S\big(\Sigma \wt y^*_\l-\wt\varphi\big)\\
&\qq+\h S\big(P \bar y^*_\l-\bar\varphi\big)-R \wt z^*_\l-\h R\bar z^*_\l\big]+\beta\Big\}\rd W(t).
\end{align*}
Comparing the diffusion terms of the above equation with those of the state equation in \rf{fbsde}, and noting that $\Sigma(\cd)\in C([0,T];\dbS^n_+)$ while $R$ and $\h R$ are positive definite, we see that $I_n+\Sigma R$ and $I_n+\Sigma\h R$ are invertible. Hence \rf{feedback-z} follows.

Finally, substituting \rf{feedback-x} and \rf{feedback-z} into the co-state equation in \rf{fbsde} and separating the centered and mean parts yields \rf{sde-y}.
This completes the proof.
\end{proof}

\br{wang923r1}
Lemmas \ref{wang923l1}--\ref{wang923l2}, together with Proposition \ref{wang922p1}, reduce the construction of the minimal norm control of Problem {\bf (NP)$_{x_0,x_T}^T$} to the determination of a parameter $\l\in\dbR^n$ satisfying \rf{wang925e8}. More precisely, the construction proceeds as follows:
\begin{enumerate}
\item [\rm(1)] Solve the Riccati equations \rf{Riccati-Sigma} and \rf{Riccati-P}, together with the MF-BSDE \rf{bsde-varphi}, to obtain $\Sigma$, $P$, $\f$, and $\beta$.

\item [\rm(2)] Determine a parameter $\l\in\dbR^n$ such that $x^*_\l(0)=x_0$.

\item [\rm(3)] For this choice of $\l$, solve the MF-SDE \rf{sde-y} to obtain $y^*_\l(\cd)$, and then compute the optimal triple $\big(x^*_\l(\cd),z^*_\l(\cd),v^*_\l(\cd)\big)$ from \rf{feedback-x}, \rf{feedback-z}, and \rf{optimality}.

\item [\rm(4)] Recover the minimal norm control $u^*(\cd)$ from \rf{wang922e8}.
\end{enumerate}
This procedure provides the basis for the numerical computation of the minimal norm control. It remains to determine a suitable $\l$ satisfying the condition in item {\rm (2)}, which will be addressed in the next subsection.
\er

\subsection{Determination of the parameter}\label{subsec-para-choice}

As indicated in Remark \ref{wang923r1}, an appropriate choice of the parameter $\l\in\dbR^n$ is crucial for ultimately solving Problem {\bf (NP)$_{x_0,x_T}^T$}. In this subsection, we provide a method to determine such a parameter.

By the feedback representation \rf{feedback-x} and the condition $x^*_\l(0)=x_0$, it follows that
\begin{align}\label{wang924e1}
x^*_\l(0)
=-P(0)\bar y^*_\l(0)-\Sigma(0)\wt y^*_\l(0)+\f(0)
=x_0.
\end{align}
By \rf{sde-y}, we have $\bar y^*_\l(0)=y^*_\l(0)=\l$, and $\wt y^*_\l(0)=0$.
Hence, \rf{wang924e1} reduces to
\begin{align}\label{wang924e2}
-P(0)\l+\f(0)=x_0.
\end{align}
Therefore, any $\l$ satisfying \rf{wang924e2} is an appropriate parameter for solving {\bf (NP)$_{x_0,x_T}^T$}.
Since $P(\cd)$ and $\f(\cd)$ are independent of $x_0$, while $x_0\in\dbR^n$ is arbitrary, the solvability of \rf{wang924e2} for every $x_0$ is equivalent to
\[
\Ran P(0)=\dbR^n.
\]
As $P(0)\in\dbR^{n\times n}$, this is equivalent to the invertibility of $P(0)$.

We next prove the invertibility of $P(0)$ under the exact controllability assumption on \rf{fmf}. To this end, we introduce an auxiliary MF-BSDE and identify its controllability Gramian with $P(0)$.

Consider the following auxiliary system:
\bel{fmf-3}
\mathrm{d}\Bx(t)=\big[\cA\wt\Bx+\h\cA \bar \Bx+\cK\wt\Bz+\h\cK\bar \Bz+\cL\wt\Bv+\h\cL\bar \Bv\big]\rd t+\Bz\rd W(t), \q t\in [0,T],
\ee
where $\Bv(\cd)\in L^2_\dbF(0,T;\dbR^{m+n})$ and
\begin{align}
&\mathcal{K}=( \BK+\Sigma S^\top)(I_n+\Sigma R)^{-1},
&&\widehat{\cK}=(\widehat{\BK}+P\widehat S^\top)(I_n+\Sigma \widehat R)^{-1},\notag\\
&\mathcal A=\BA+\Sigma Q-\mathcal{K}\Sigma S,
&&\widehat{\mathcal A}=\widehat{\mathbf A}+P\widehat Q-
\widehat{\mathcal K}\Sigma \widehat S,\notag\\
&\cL=\big(\BL,\, (-\Sigma,\cK\Sigma)\Theta^{1/2}\big),
&&\h\cL=\big(\h\BL,\, (-P,\h\cK\Sigma)\h\Theta^{1/2}\big).\label{auxi-matrix}
\end{align}
Here $\Theta^{1/2}$ and $\h\Theta^{1/2}$ denote the symmetric positive semidefinite square roots of the matrices $\Theta$ and $\h\Theta$ defined in \rf{wang923e1}, respectively. Since $r=m-n$, both $\cL$ and $\h\cL$ belong to $\dbR^{n\times(m+n)}$.

We first establish the exact controllability of the auxiliary system \rf{fmf-3} and recall its Gramian characterization.

\bl{wang924l1}
Suppose that \rf{fmf} is exactly controllable on $[0,T]$. Then the auxiliary system \rf{fmf-3} is also exactly controllable on $[0,T]$. Furthermore, the exact controllability of \rf{fmf-3} on $[0,T]$ is equivalent to the invertibility of the following Gramian:
\begin{align}\label{wang924e3}
G_{\tt aux}\deq\me\int_0^T\big[
\wt\Phi\cL\cL^\top\wt\Phi^\top
+\bar\Phi\h{\cL}\h{\cL}^\top\bar\Phi^\top
\big]\rd t,
\end{align}
where $\Phi(\cd)$ is the fundamental matrix associated with \rf{fmf-3}, satisfying
\begin{equation}\label{sde-Phi}
\left\{
\begin{aligned}
&\mathrm{d}\Phi(t)=
-\big[\wt\Phi\cA+\bar \Phi\h{\cA}\big]\rd t
-\big[\wt\Phi\cK+\bar \Phi\h{\cK}\big]\rd W(t),\q t\in [0,T],\\
&\Phi(0)=I_n.
\end{aligned}
\right.
\end{equation}
\el

\begin{proof}
Since \rf{fmf} is exactly controllable on $[0,T]$, it is also terminally controllable on $[0,T]$. Hence, Theorem \ref{rank-condition} and Proposition \ref{wang922p1} imply that \rf{fmf-2} is exactly controllable on $[0,T]$.

Fix any $x_0\in\dbR^n$ and $x_T\in L^2_{\mf_T}(\Omega;\dbR^n)$. 
Then there exists a control $v(\cd)\in L^2_\dbF(0,T;\dbR^r)$ such that the corresponding solution $\big(x(\cd),z(\cd)\big)$ of \rf{fmf-2}, with the terminal condition $x(T)=x_T$, satisfies $x(0)=x_0$.
Define
\begin{align*}
\Bv=
\begin{pmatrix}
v\\
\Theta^{1/2}\begin{pmatrix}\wt x\\\wt z\end{pmatrix}
+\h\Theta^{1/2}\begin{pmatrix}\bar x\\\bar z\end{pmatrix}
\end{pmatrix}.
\end{align*}
Then $\Bv(\cd)\in L^2_\dbF(0,T;\dbR^{m+n})$ and
\begin{align*}
\wt\Bv=
\begin{pmatrix}
\wt v\\
\Theta^{1/2}\begin{pmatrix}\wt x\\ \wt z\end{pmatrix}
\end{pmatrix},
\qquad
\bar\Bv=
\begin{pmatrix}
\bar v\\
\h\Theta^{1/2}\begin{pmatrix}\bar x\\ \bar z\end{pmatrix}
\end{pmatrix}.
\end{align*}
Using \rf{auxi-matrix}, together with $\cK(I_n+\Sigma R)=\BK+\Sigma S^\top$ and $\h\cK(I_n+\Sigma\h R)=\h\BK+P\h S^\top$, a direct calculation gives
\begin{align*}
\cA\wt x+\h\cA\bar x+\cK\wt z+\h\cK\bar z
+\cL\wt\Bv+\h\cL\bar\Bv
=
\BA\wt x+\h\BA\bar x+\BK\wt z+\h\BK\bar z
+\BL\wt v+\h\BL\bar v.
\end{align*}
Thus $\big(x(\cd),z(\cd)\big)$ also satisfies \rf{fmf-3} with the terminal condition $x(T)=x_T$ and the control $\Bv(\cd)$. By the uniqueness of the solution to the corresponding MF-BSDE, we have
\[
\big(\Bx(\cd),\Bz(\cd)\big)=\big(x(\cd),z(\cd)\big),
\]
and consequently $\Bx(0)=x(0)=x_0$.
Since $x_0$ and $x_T$ are arbitrary, \rf{fmf-3} is exactly controllable on $[0,T]$.

Finally, the equivalence between the exact controllability of \rf{fmf-3} and the invertibility of the Gramian \rf{wang924e3} follows from \cite[Theorem 4.1]{Yu21}.
This completes the proof.
\end{proof}

We next identify the Gramian $G_{\tt aux}$ with the value $P(0)$ of the Riccati solution \rf{Riccati-P}.

\bt{inver-P}
Let $P(\cd)$ be the solution to the Riccati equation \rf{Riccati-P} and $G_{\tt aux}$ be the controllability Gramian \rf{wang924e3} of the auxiliary system \rf{fmf-3}. Then
\bel{P=G}
P(0)=G_{\tt aux}.
\ee
Consequently, if \rf{fmf} is exactly controllable on $[0,T]$, then $P(0)$ is invertible.
\et

\begin{proof}
We divide the proof into three steps.

\no{\bf Step 1.}
Taking expectations in \rf{sde-Phi}, we have
\begin{equation*}
\left\{
\begin{aligned}
&\dot{\bar\Phi}=-\bar \Phi\h{\cA},\qq \mbox{on } [0,T],\\
&\bar\Phi(0)=I_n.
\end{aligned}
\right.
\end{equation*}
Applying the fundamental theorem of calculus to $\bar\Phi(\cd)(P(\cd)-\Sigma(\cd))\bar\Phi(\cd)^\top$ and using
$P(T)=\Sigma(T)=0$, $\bar \Phi(0)=I_n$, we obtain
\begin{equation}\label{wang924e4}
\begin{aligned}
P(0)-\Sigma(0)=\int_0^T \big[\bar\Phi\widehat{\mathcal A}(P-\Sigma )\bar\Phi^\top
+\bar\Phi(P-\Sigma )\widehat{\mathcal A}^\top\bar\Phi^\top
-\bar\Phi(\dot P-\dot \Sigma )\bar\Phi^\top\big] \rd t.
\end{aligned}
\end{equation}

On the other hand, applying It\^o's formula to $\Phi(\cd)\Sigma(\cd)\Phi^\top(\cd)$ gives
\begin{align*}
\mathrm{d}(\Phi \Sigma \Phi^\top)
&=-\Big[\big(\widetilde\Phi\mathcal A+\bar\Phi\widehat{\mathcal A}\big)\Sigma \Phi^\top +\Phi \Sigma \big(\widetilde\Phi\mathcal A+\bar\Phi\widehat{\mathcal A}\big)^\top
-\Phi\dot \Sigma \Phi^\top \\
&\qq\q-\big(\widetilde\Phi\mathcal{K}+\bar\Phi\widehat{\mathcal K}\big)
\Sigma \big(\widetilde\Phi\mathcal{K}+\bar\Phi\widehat{\mathcal K}\big)^\top\Big] \rd t
+\cds\rd W(t).
\end{align*}
Since $\Sigma(T)=0$ and $\Phi(0)=I_n$, taking expectations gives
\bel{wang924e5}
\bal
\Sigma(0)&=\me\int_0^T\Big\{\wt\Phi\big[ \cA\Sigma+\Sigma\cA^\top-\cK\Sigma\cK^\top-\dot\Sigma\big]\wt\Phi^\top\\
&\qq\q+\bar\Phi\big[ \h\cA\Sigma+\Sigma\h\cA^\top-\h\cK\Sigma\h\cK^\top-\dot\Sigma\big]\bar\Phi^\top\Big\}\rd t.
\eal
\ee
Here we use the fact that $\Phi=\widetilde\Phi+\bar\Phi$, $\me[\widetilde\Phi]=0$, and $\me\big[\wt \Phi M\bar \Phi^\top\big]=\me\big[\bar \Phi M\wt \Phi^\top\big]=0$ for any $M\in\dbR^{n\times n}$.
Consequently, adding \rf{wang924e4} and \rf{wang924e5} yields
\bel{wang924e6}
\bal
P(0)&=\me\int_0^T\Big\{\wt\Phi\big[ \cA\Sigma+\Sigma\cA^\top-\cK\Sigma\cK^\top-\dot\Sigma\big]\wt\Phi^\top\\
&\qq+\bar\Phi\big[ \h\cA P+P \h\cA^\top-\h\cK \Sigma \h\cK^\top-\dot P\big]\bar\Phi^\top\Big\}\rd t.
\eal
\ee

\no{\bf Step 2.} We claim that
\bel{wang924e7}
\bal
&\dot \Sigma -\mathcal A\Sigma -\Sigma \mathcal A^\top+
\mathcal{K}\Sigma \mathcal{K}^\top=-\mathcal{L}\mathcal{L}^\top,\\
\mbox{and}\q&
\dot P-
\widehat{\mathcal A}P-
P\widehat{\mathcal A}^\top+
\widehat{\mathcal K}\Sigma \widehat{\mathcal K}^\top
=-\widehat{\mathcal{L}}\widehat{\mathcal{L}}^\top.
\eal
\ee
We only prove the first identity, since the second one follows in the same vein.
By \rf{auxi-matrix},
\begin{align*}
&\dot \Sigma -\mathcal A\Sigma -\Sigma \mathcal A^\top+
\mathcal{K}\Sigma \mathcal{K}^\top \notag\\
&\qq=\dot \Sigma -\mathbf A\Sigma -\Sigma \mathbf A^\top-2\Sigma Q\Sigma 
+\mathcal{K}\Sigma S\Sigma +\Sigma S^\top \Sigma \mathcal{K}^\top
+\mathcal{K}\Sigma \mathcal{K}^\top.
\end{align*}
Since $ \BK+\Sigma S^\top=\mathcal{K}(I+\Sigma R)$, the Riccati equation \rf{Riccati-Sigma} can be rewritten as 
\[
\dot \Sigma -\mathbf A\Sigma -\Sigma \mathbf A^\top-\Sigma Q\Sigma +
 \BL\BL^\top
+\mathcal{K}\Sigma \mathcal{K}^\top
+\mathcal{K}\Sigma R\Sigma \mathcal{K}^\top=0.
\]
Hence, comparing the above two equations gives
\begin{align}
&\dot \Sigma -\mathcal A\Sigma -\Sigma \mathcal A^\top+
\mathcal{K}\Sigma \mathcal{K}^\top\notag\\
&\qq=-\Sigma Q\Sigma - \BL\BL^\top
-\mathcal{K}\Sigma R\Sigma \mathcal{K}^\top
+\mathcal{K}\Sigma S\Sigma +\Sigma S^\top \Sigma \mathcal{K}^\top.\label{wang924e10}
\end{align}

On the other hand, by \eqref{wang923e1}, it holds that
\begin{align*}
(-\Sigma,\mathcal{K}\Sigma )\Theta(-\Sigma,\mathcal{K}\Sigma )^\top
=\Sigma Q\Sigma+\mathcal{K}\Sigma R\Sigma \mathcal{K}^\top -\mathcal{K}\Sigma S\Sigma -\Sigma S^\top \Sigma \mathcal{K}^\top.
\end{align*}
Therefore, by the definition of $\cL$ in \rf{auxi-matrix},
\begin{align}
\mathcal{L}\mathcal{L}^\top
=\Sigma Q\Sigma +\BL\BL^\top
+\mathcal{K}\Sigma R\Sigma \mathcal{K}^\top-\mathcal{K}\Sigma S\Sigma -\Sigma S^\top \Sigma \mathcal{K}^\top.
\label{wang924e11}
\end{align}
Now the assertion on $\Sigma$ in \rf{wang924e7} can be established by combining \rf{wang924e10} and \rf{wang924e11}.
The $P$-part in \rf{wang924e7} follows analogously.

\no
{\bf Step 3.} A combination of \rf{wang924e6}, \rf{wang924e7} and the definition of the Gramian $G_{\tt aux}$ in \rf{wang924e3} implies that
\begin{align*}
P(0)=\me\int_0^T\big[\wt\Phi \cL\cL^\top \wt\Phi^\top
+\bar\Phi\h\cL\h\cL^\top\bar\Phi^\top\big]\rd t=G_{\tt aux}.
\end{align*}
That settles the assertion \rf{P=G}.

Under the exact controllability assumption on the original system \rf{fmf}, Lemma \ref{wang924l1} implies that the auxiliary system \rf{fmf-3} is exactly controllable. Hence its Gramian $G_{\tt aux}$ is invertible. Together with \rf{P=G}, this yields the invertibility of $P(0)$.
This completes the proof.
\end{proof}


Under the exact controllability assumption on \rf{fmf}, Theorem \ref{inver-P} implies that the parameter $\l$ satisfying \rf{wang924e2} is uniquely determined by
\bel{wang924e12}
\l=P(0)^{-1} \big(\f(0)-x_0\big).
\ee
We can therefore complete the construction described in Remark \ref{wang923r1} and obtain the following explicit representation of the minimal norm control of Problem {\bf (NP)$_{x_0,x_T}^T$}.

\bc{wang924c1}
Suppose that the MF-SDE system \rf{fmf} is exactly controllable on $[0,T]$. Let $\Sigma(\cd), P(\cd)$, and $\big(\f(\cd),\b(\cd)\big)$ be the solutions to the Riccati equations \rf{Riccati-Sigma}, \rf{Riccati-P}, and the MF-BSDE \rf{bsde-varphi}, respectively. Define $\l$ by \rf{wang924e12}, and let $y^*_\l(\cd)$ be the solution of \rf{sde-y} corresponding to this choice of $\l$. Then the optimal triple $\big(x^*_\l(\cd),z^*_\l(\cd),v^*_\l(\cd)\big)$ of Problem {\bf (OP)$_{x_0,x_T,\l}^T$} is given by
\bel{OP-triple}
\left\{
\bal
&x^*_\l=-P\bar y^*_\l-\Sigma \wt y^*_\l+\f,\\
&z^*_\l=(I+\Sigma R)^{-1}\big[\Sigma  \big(\BK^\top+S\Sigma\big) \wt y^*_\l-\Sigma S\wt \varphi+\wt\beta\big]\\
&\qq\q+(I+\Sigma \h R)^{-1}\big[\Sigma  \big(\h \BK^\top+\h S P\big) \bar y^*_\l-\Sigma \h S\bar \varphi+\bar\beta\big], \\
&v^*_\l=\BL^\top \wt y^*_\l+\h\BL^\top \bar y^*_\l.
\eal
\right.
\ee
Moreover, the minimal norm control $u^*(\cd)$ of Problem {\bf (NP)$_{x_0,x_T}^T$} is given by
\bel{NP-control}
u^*=D^\top(DD^\top)^{-1} \big[\wt z^*_\l-C\wt x^*_\l\big]+U \wt v^*_\l 
+ \widehat D^\top(\widehat D\widehat D^\top)^{-1} \big[\bar z^*_\l-\h C \bar x^*_\l\big]+\h U \bar v^*_\l,
\ee
where $U$ and $\h U$ are the matrices fixed in \rf{decom}.
\ec

The preceding results also yield the following characterization of exact controllability.

\bp{wang929p1}
The MF-SDE system \rf{fmf} is exactly controllable on $[0,T]$ if and only if
\[
\Rank(D)=\Rank(D+\bar D)=n
\]
and the solution $P(\cd)$ of the Riccati equation \rf{Riccati-P}, defined under this rank condition, satisfies 
\[
P(0)>0.
\]
\ep

\begin{proof}
Suppose first that the MF-SDE system \rf{fmf} is exactly controllable on $[0,T]$. Then it is also terminally controllable on $[0,T]$. Hence, by Theorem \ref{rank-condition}, $\Rank(D)=\Rank(D+\bar D)=n$. Under this rank condition, Lemma \ref{wang924l1} implies that the auxiliary system \rf{fmf-3} is exactly controllable on $[0,T]$, and hence its controllability Gramian $G_{\tt aux}$ defined in \rf{wang924e3} is positive definite. Since Theorem \ref{inver-P} gives $P(0)=G_{\tt aux}$, it follows that $P(0)$ is positive definite.

Conversely, suppose that $\Rank(D)=\Rank(D+\bar D)=n$ and that $P(0)$ is positive definite. Fix any $x_0\in\dbR^n$ and $x_T\in L^2_{\mf_T}(\Omega;\dbR^n)$, and let $\big(\f(\cd),\beta(\cd)\big)$ be the solution of \rf{bsde-varphi} with terminal condition $\f(T)=x_T$. Set $\l=P(0)^{-1}\big(\f(0)-x_0\big)$.
By Lemmas \ref{wang923l1} and \ref{wang923l2}, the corresponding optimal triple $\big(x^*_\l(\cd),z^*_\l(\cd),v^*_\l(\cd)\big)$
of Problem {\bf (OP)$_{x_0,x_T,\l}^T$} satisfies the MF-BSDE \rf{fmf-2}, with
\[
x^*_\l(0)
=-P(0)\l+\f(0)
=x_0,
\qquad
x^*_\l(T)=x_T.
\]
Since $x_0$ and $x_T$ are arbitrary, the MF-BSDE system \rf{fmf-2} is exactly controllable on $[0,T]$. Proposition \ref{wang922p1} then implies that the MF-SDE system \rf{fmf} is exactly controllable on $[0,T]$. This completes the proof.
\end{proof}

\br{wang929r1}
Proposition \ref{wang929p1} provides a deterministic Riccati-based criterion for the exact controllability of the constant-coefficient MF-SDE system \rf{fmf}. In contrast with the Gramian characterization in Lemma \ref{wang924l1}, which requires solving the matrix-valued stochastic equation \rf{sde-Phi} and evaluating the expectation in \rf{wang924e3}, the condition
$P(0)>0$
can be checked by solving the two deterministic Riccati equations \rf{Riccati-Sigma}--\rf{Riccati-P} backward in time. Moreover, these Riccati equations are already required in the construction of the minimal norm control in Corollary \ref{wang924c1}. Thus, the controllability test introduces essentially no additional computational procedure beyond that used to construct the optimal control.

This characterization is also useful for numerical verification of exact controllability. Let $P_h(0)$ be a numerical approximation of $P(0)$, and suppose that an error bound
\[
\|P_h(0)-P(0)\|_2\leq\varepsilon_h
\]
is available. Here $\|\cdot\|_2$ denotes the matrix spectral norm induced by the Euclidean norm on $\dbR^n$. Since $P(0)$ is symmetric, Weyl's eigenvalue perturbation inequality (see, e.g., \cite[Theorem 4.3.1]{Horn-Johnson90}) gives
\[
\lambda_{\min}\big(P(0)\big)
\geq
\lambda_{\min}\Big(\frac{P_h(0)+P_h(0)^\top}{2}\Big)-\varepsilon_h.
\]
Consequently, under the rank condition $\Rank(D)=\Rank(D+\bar D)=n$, the computable condition
\[
\lambda_{\min}\Big(\frac{P_h(0)+P_h(0)^\top}{2}\Big)>\varepsilon_h
\]
guarantees the exact controllability of \rf{fmf}. 
\er

%
%
%

Finally, the Gramian identity \rf{P=G} also yields an explicit minimal norm control for the auxiliary system \rf{fmf-3}.

\bc{wang925c1}
Suppose that the MF-SDE system \rf{fmf} is exactly controllable on $[0,T]$, and let $\l$ be given by \rf{wang924e12}. Then the unique minimal norm control of the auxiliary system \rf{fmf-3} satisfying $\Bx(0)=x_0$ and $\Bx(T)=x_T$ is given by
\bel{wang925e1}
\Bv^*
=\big[
\cL^\top \wt\Phi^\top
+\h\cL^\top\bar\Phi^\top
\big]\l.
\ee
More precisely,
\begin{align*}
\|\Bv^*(\cd)\|_{L^2_\dbF(0,T;\dbR^{m+n})}
&=\inf\Big\{\|\Bv(\cd)\|_{L^2_\dbF(0,T;\dbR^{m+n})}\,\big|\, \Bv(\cd)\in L^2_\dbF(0,T;\dbR^{m+n}),\\
&\qq\qq\mbox{ and } \Bx(0;x_T,\Bv(\cd))=x_0\Big\}.
\end{align*}
Moreover,
\bel{wang925e2}
\|\Bv^*(\cd)\|_{L^2_\dbF(0,T;\dbR^{m+n})}^2
=
\big\langle
P(0)^{-1}\big(\f(0)-x_0\big),
\f(0)-x_0
\big\rangle
=
\langle P(0)\l,\l\rangle.
\ee
\ec

\begin{proof}
By \rf{bsde-varphi} and the definitions in \rf{auxi-matrix}, $\big(\f(\cd),\b(\cd)\big)$ is precisely the solution pair of the auxiliary system \rf{fmf-3} corresponding to the terminal condition $x_T$ and the zero control. Hence,
\[
\Bx(0;x_T,0)=\f(0).
\]
Applying the minimal norm controllability result in \cite[Theorem 4.1 and Corollary 4.2]{Chen-Yu24} to \rf{fmf-3} gives
\[
\Bv^*
=-\big[
\cL^\top\wt\Phi^\top
+\h\cL^\top\bar\Phi^\top
\big]
G_{\tt aux}^{-1}\big(x_0-\f(0)\big).
\]
Using $G_{\tt aux}=P(0)$ from Theorem \ref{inver-P} and \rf{wang924e12} yields \rf{wang925e1}. The minimal norm identity \rf{wang925e2} follows immediately.
This completes the proof.
\end{proof}

We conclude this section with an explicit example illustrating the Riccati-based representation obtained above. 

\bex{explicit-example}
Let $T=1$, $n=1$, $m=2$, $x_0=0$, $x_T=1$, and set
\begin{align*}
A=-\frac12,\,
\bar A=\frac12,\,
B=(1,-1),\,
\bar B=(-1,2),\,C=1,\,
\bar C=-1,\,
D=(0,1),\,
\bar D=(0,0).
\end{align*}
Then $\Rank D=\Rank\h D=1$, and we may choose $U=\h U=(1,\,0)^\top$.

For these coefficients, the solutions of the Riccati equations \rf{Riccati-Sigma} and \rf{Riccati-P} are $\Sigma(t)=1-t$ and $P(t)=(1-t)-\ln(2-t)$ with $t\in [0,1]$, respectively. In particular,
\[
P(0)=1-\ln2>0.
\]
Hence, by Proposition \ref{wang929p1}, the MF-SDE system \rf{fmf} is exactly controllable on $[0,1]$.
Then the solution of the MF-BSDE \rf{bsde-varphi} is $\big(\f(t),\beta(t)\big)=(1,0)$ for $t\in[0,1]$.
Therefore, by \rf{wang924e12},
\[
\l
=P(0)^{-1}\big(\f(0)-x_0\big)
=\frac{1}{1-\ln2}.
\]
For this choice of $\l$, the MF-SDE \rf{sde-y} becomes
\beq
\left\{
\bal
&\mathrm{d} y^*_\l(t)
=-\frac12\wt y^*_\l(t)\rd t
+\Big[\wt y^*_\l(t)-\frac{\bar y^*_\l(t)}{2-t}\Big]\rd W(t),
\qquad t\in[0,1],\\
&y^*_\l(0)=\l.
\eal
\right.
\eeq
Taking expectations gives $\bar y^*_\l(t)=\l$. Consequently, the centered process $\wt y^*_\l(\cd)$ satisfies a linear SDE, whose explicit solution is
\begin{align*}
\wt y^*_\l(t)
=\l e^{W(t)-t}
\Big[\int_0^t\frac{e^{-W(s)+s}}{2-s}\rd s
-\int_0^t\frac{e^{-W(s)+s}}{2-s}\rd W(s)
\Big].
\end{align*}

Using \rf{optimality}, \rf{feedback-z}, and \rf{NP-control}, we obtain the unique minimal norm control
\begin{align*}
u^*(t)
=\Big(\wt y^*_\l(t),\, \frac{\l(1-t)}{2-t}\Big)^\top,\qq t\in[0,1].
\end{align*}
Moreover, It\^o's formula applied to $|\wt y^*_\l(t)|^2$ gives
$\me\big[|\wt y^*_\l(t)|^2\big]
=\l^2\big(\frac{1}{2-t}-\frac12\big)>0$,
for $t\in(0,1]$.
Thus, although both $x_0=0$ and $x_T=1$ are deterministic, the first component of the minimal norm control is genuinely stochastic.
\ex

\br{wang929r2}
The above example also illustrates the need for developing numerical methods for NP problems of MF-SDE systems. Even when both the initial state and the terminal target are deterministic, the minimal norm control may contain a nontrivial stochastic component and generally does not admit a simple closed-form representation. On the other hand, the constructive approach developed in this paper expresses the minimal norm control in terms of Riccati equations, an MF-BSDE, and a decoupled MF-SDE, thereby providing a feasible framework for its numerical computation. The design and convergence analysis of the corresponding numerical algorithms will be investigated in our future work.
\er



\end{document}